\documentclass[12pt]{amsart}
\usepackage{amscd,amsmath,amssymb,amsthm,enumerate,times,ifthen}
\usepackage{epsfig,epic,eepic}
\usepackage{graphicx}
\usepackage[utf8]{inputenc}
\usepackage[T1]{fontenc}
\usepackage{tikz}
\usetikzlibrary{matrix,arrows}
\usepackage[mathscr]{eucal}
\newtheorem{theorem}{Theorem}
\newtheorem{Theorem}{Theorem}

\newtheorem{Lemma}{Lemma}

\newtheorem{Prposition}{Proposition}

\newtheorem{Corollary}{Corollary}

\newtheorem*{theoremn}{Theorem}
\def\Thmn#1{\begin{theoremn}#1\end{theoremn}}

\newtheorem*{lemman}{Lemma}
\def\Lemman#1{\begin{lemman}#1\end{lemman}}

\newtheorem*{corn}{Corollary}
\def\Corn#1{\begin{corn}#1\end{corn}}

\theoremstyle{definition}
\newtheorem*{Def}{Definition}

\newcommand{\N}{\mathbb N}

\newcommand{\R}{\mathbb R}

\newcommand{\eps}{{\varepsilon}}
\newcommand{\id}{\mathop{\hbox{\rm id}}}

\newcommand{\diam}{\mathop{\hbox{\rm diam}}}
\newcommand{\orbit}{\mathcal{O}}
\renewcommand{\phi}{\varphi}
\newcommand{\eq}[1]{\eq{#1}}
\newcommand{\Eq}[2]{\ifthenelse{\equal{#1}{*}}
{\begin{equation*}\begin{aligned}#2\end{aligned}\end{equation*}}
{\begin{equation}\label{#1}\begin{aligned}#2\end{aligned}\end{equation}}}

\begin{document}

\date{\today}

\begin{flushright}
\emph{Publ. Math. Debrecen (Submitted)}
\end{flushright} 

\title{The contraction principle in extended context}
\dedicatory{Dedicated to the 60th birthday of Professor Zsolt P\'ales}
\author[M. Bessenyei]{Mih\'aly Bessenyei}
\address{Institute of Mathematics, University of Debrecen, H-4010 Debrecen, Pf.\ 12, Hungary}
\email{besse@science.unideb.hu}



\thanks{This research has been supported by the Hungarian Scientific Research Fund (OTKA) Grants K--111651.}

\begin{abstract}
There are several extensions of the classical Banach Fixed Point Theorem in technical literature. A branch of
generalizations replaces usual contractivity by weaker but still effective assumptions. Our note follows this
stream, presenting an elementary proof for a known fixed point result. Some applications are also considered.
\end{abstract}

\maketitle

\section{Introduction}

Although the Contraction Principle appears partly in the method of successive approximation trough the works of
Cauchy \cite{Cau1835}, Liouville \cite{Lio1837}, and Picard \cite{Pic1890}, its abstract and powerful version
is due to Banach \cite{Ban22} and Cacciopoli \cite{Cac30}. This form of the Contraction Principle, quoted as
the Banach Fixed Point Theorem, states that any contraction of a complete metric space has exactly one fixed
point.

We can say that the Contraction Principle has made an explosion in contemporary mathematics. It has become the
corner stone of Iterative Fixed Point Theory, has initiated new branches of important generalizations, and has
found its legacy in many fields of mathematics via applications. The monographs by Berinde \cite{Ber07}, by
Granas and Dugundji \cite{GraDug03}, by Rus, Petru\c{s}el and Petru\c{s}el \cite{RusPetPet08}, and by Zeidler
\cite{Zei86} give an excellent and detailed demonstration of this fact.

An important stream of the extensions of the Contraction Principle replaces contractivity by weaker but still
effective properties. A rich overview of such generalizations can be found in the survey of Rhoades \cite{Rho77}.
Among them, let us recall the pioneer works of Boyd and Wong \cite{BoyWon69}, of Browder \cite{Brow68}, and of
Matkowski \cite{Mat75}. These investigations present fixed point results for \emph{nonlinear contractions}. An
other direction concerns \emph{linear quasicontractions}, and was studied first by \'Ciri\'c \cite{Cir74} and
also by Kannan \cite{Kan69}. For precise details, the Reader can have a quick look at the corollaries of this note.
Surprisingly, applying \emph{(nonlinear) quasicontractions}, the above mentioned cases get a common point of view:
Independently and simultaneously, Heged\H us and Szil\'agyi \cite{HegSzi80} and Walter \cite{Wal81} have given
fixed point results for this general setting.

These extensions have not merely theoretical importance, but enjoy many didactic and philosophical aspect, as
well. For those ones, who read through and understand the various ideas, open a colorful perspective of the
Contraction Principle. 

The aim of the present note is to re-discover the theorem of Heged\H us, Szil\'agyi and Walter, this wonderful and
hidden treasure of Fixed Point Theory. Requiring a reasonable extra assumption on comparison functions, an elementary
and self-contained approach can be given. This approach does not exceed the standard tools of classical analysis and
the basic facts on metric spaces, and hence, according to our hope, might have some impact even in education.

\section{Notations, notions, basic facts}

Throughout this note, $\R_+$ and $\N$ stand for the set of nonnegative reals and natural numbers, respectively.
If $f$ is a self-map of a nonempty set, then the \emph{composite iterates} are defined via the usual recursion
$f^{n+1}=f\circ f^n$ under the convention $f^0=\id$. The \emph{orbit} and the \emph{double orbit} induced by $f$
are defined in the next way:
\Eq{*}{
 \orbit(x):=\{f^n(x)\mid n\in\N\cup\{0\}\};\qquad
  \orbit(x,y):=\orbit(x)\cup\orbit(y).}

Our investigations focus on nonlinear quasicontractions, when the distance of images is controlled by a given
function of double orbits induced by the original points. 

\begin{Def}
Under a \emph{comparison function} we mean an increasing, upper semicontinuous function $\phi\colon\R_+\to\R_+$
which fulfills the properties $\phi(0)=0$ and $\phi(t)<t$ for $t>0$. Let $(X,d)$ be a metric space. A mapping
$T\colon X\to X$ is called a \emph{weak quasicontraction} with comparison function $\phi$ (or briefly: a
\emph{weak $\phi$-quasicontraction}) if it induces bounded orbits and, for all $x,y\in X$,
\Eq{*}{
 d(Tx,Ty)\le\phi\bigl(\diam\orbit(x,y)\bigr).}
Similarly, under a \emph{strong $\phi$-quasicontraction} we mean a mapping $T\colon X\to X$ fulfilling the property
below for all $x,y\in X$:
\Eq{*}{
 d(Tx,Ty)\le\phi\bigl(\diam\{x,y,Tx,Ty\}\bigr).}
\end{Def}

Clearly, if a strong quasicontraction induces bounded orbits, then it is a weak quasicontraction. However, there
exist strong quasicontractions that induce unbounded orbits \cite{Bes15a}, and hence they are not considered weak
quasicontractions in our framework. Of course, the definition above subsumes the notions of linear quasicontractions
and nonlinear contractions mentioned in the Introduction. The most unusual feature of a weak/strong quasicontraction
is that its continuity cannot be guaranteed, even in the most simple cases \cite{Kan69}. This phenomenon makes
impossible to apply the standard methods: New ideas have to be developed for establishing fixed point theorems. 

Note that \emph{a weak quasicontraction can have at most one fixed point}. Indeed, assume that $x_0$ and $y_0$ are
distinct fixed points of a weak $\phi$-quasicontraction $T$. Then, $d(x_0,y_0)>0$ and hence
\Eq{*}{
 d(x_0,y_0)=d(Tx_0,Ty_0)\le\phi\bigl(\diam\orbit(x_0,y_0)\bigr)=\phi\bigl(d(x_0,y_0)\bigr)<d(x_0,y_0).}
A similar argument shows that \emph{the fixed point of a strong quasicontraction is unique} provided that it exists
at all.

In the sequel, let us present two less obvious properties for comparison functions and weak quasicontractions. The
first one is well-known (see for example \cite{Bes15a}). For the Readers' convenience, we sketch here the proof.  

\Lemman{The sequence of composite iterates of a comparison function tends to zero pointwise.}

\begin{proof}
If $\phi$ is a comparison function, then the limit property obviously holds at $t=0$. Fix $t>0$. If $\phi(t)=0$,
then $\phi^n(t)=0$ for all $n\in\N$; if $\phi(t)\neq 0$, then $\phi^2(t)<\phi(t)$ follows after iterating the
inequality $\phi(t)<t$. Proceeding by induction, we get that $(\phi^n)$ is decreasing on $\R_+$. On the other
hand, it is bounded from below. Therefore, $f=\lim_{n\to\infty}\phi^n$ exists and takes nonnegative values.
Assume indirectly that $f(t)>0$ for some $t>0$. Then, by upper semicontinuity, we arrive at
\Eq{*}{
 f(t)=\lim_{n\to\infty}\phi^{n+1}(t)=
  \limsup_{n\to\infty}\phi\bigl(\phi^n(t)\bigr)\le
   \phi\bigl(\limsup_{n\to\infty}\phi^n(t)\bigr)=
    \phi\bigl(f(t)\bigr)<f(t).}
This contradiction completes the proof.
\end{proof}

Observe that monotonicity has played no role here. However, the proof of the next auxiliary lemma enlightens the
importance of this property.

\Lemman{If $T$ is a weak $\phi$-quasicontraction, then $T^n$ is a weak $\phi^n$-quasicontraction.}

\begin{proof}
Fix $k,l\in\N$ and let $T^kx,T^ly\in\orbit(x,y)$. Then the contractivity of $T$ and the monotonicity of $\phi$
gives
\Eq{*}{
 d(T^kx,T^ly)=d(TT^{k-1}x,TT^{l-1}y)
  \le\phi\bigl(\diam\orbit(T^{k-1}x,T^{l-1}y)\bigr)
   \le\phi\bigl(\diam\orbit(x,y)\bigr).}
A similar argument leads to
\Eq{*}{
 d(T^kx,T^lx)\le\phi\bigl(\diam\orbit(T^{k-1}x,T^{l-1}x)\bigr)
  \le\phi\bigl(\diam\orbit(x)\bigr)
   \le\phi\bigl(\diam\orbit(x,y)\bigr).}
For $T^ky,T^ly\in\orbit(x,y)$, the same upper estimation can be obtained in the same way. That is, we have the next
inequality between the double orbits of $\{x,y\}$ and $\{Tx,Ty\}$: 
\Eq{*}{
 \diam\orbit(Tx,Ty)=
  \sup_{k,l\in\N}\{d(T^kx,T^ly),d(T^kx,T^lx),d(T^ky,T^ly)\}
   \le\phi\bigl(\diam\orbit(x,y)\bigr).}

To complete the proof, we apply induction. For $n=1$, the statement holds trivially. Assume that it also remains
true for some $n\in\N$. Then, applying the inequality above and the monotonicity of the comparison function,
\Eq{*}{
 d(T^{n+1}x,T^{n+1}y)=d(T^nTx,T^nTy)\le\phi^n\bigl(\diam\orbit(Tx,Ty)\bigr)
  \le\phi^{n+1}\bigl(\diam\orbit(x,y)\bigr)}
follows, which was to be proved.
\end{proof}

Note also that an analogous statement cannot be formulated for strong quasicontractions. The reason for this, in
particular, is that the iterates of a strong quasicontraction does not generate necessarily bounded orbits, even if
the original mapping does.

\section{The main result}

The main result of this note presents a fixed point theorem for weak quasicontractions of complete metric spaces.
The proof has three stages. The first and standard one is devoted to show the Cauchy property of sequence of iterates.
Then completeness gives the existence of a limit point. The forthcomings bring novelty: In the second step we show,
that the iterates on the limit point tend also to this limit. In fact, it is an immediate consequence in the standard
cases but now, in lack of continuity, will have a particular importance. The final step concludes that the orbit induced
by the limit has zero diameter, which is the desired fixed point property. The technical lemmas of the previous section
play a key role in the arguments.  

\Thmn{Any weak quasicontraction of a complete metric space has a unique fixed point. Moreover, the sequence of iterates
at any point converges to this fixed point.}

\begin{proof}
Let $(X,d)$ complete metric space and let $T\colon X\to X$ be a weak $\phi$-quasicontraction. Fix $x\in X$ arbitrarily.
We show that $(T^nx)_{n=0}^{\infty}$ has the Cauchy property. The boundedness of orbits and the pointwise convergence
property of $(\phi^n)$ provides that, for all $\eps>0$, there exists $n_0\in\N$ such that 
$\phi^{n_0}\bigl(\diam\orbit(x)\bigr)<\eps/2$. If $n>n_0$, then
\Eq{*}{
 d(T^{n_0}x,T^nx)\le\phi^{n_0}\bigl(\diam\orbit(x,T^{n-n_0}x)\bigr)=
  \phi^{n_0}\bigl(\diam\orbit(x)\bigr)<\eps/2.}
Hence the use of the triangle inequality leads to $d(x_n,x_m)<\eps$ for all $n,m>n_0$. In other words, the iterates
of $x$ defines a Cauchy sequence.

The completeness guarantees the existence of $x_0\in X$ such that $T^nx\to x_0$. Our claim is that $T^nx_0\to x_0$
also holds. Indeed,
\Eq{*}{
 d(x_0,T^nx_0)\le d(x_0,T^nx)+d(T^nx,T^nx_0)
  \le d(x_0,T^nx)+\phi^n\bigl(\diam\orbit(x,x_0)\bigr)
   \to 0}
as $n\to\infty$.

To complete the proof, we suffice to show that $\diam\orbit(x_0)=0$. Suppose to the contrary that this is not the
case. Then, for all $n,k\in\N$, we arrive at
\Eq{*}{
 d(T^nx_0,T^{n+k}x_0)\le\phi^n\bigl(\diam\orbit(x_0,T^kx_0)\bigr)
  =\phi^n\bigl(\diam\orbit(x_0)\bigr)
   \le\phi\bigl(\diam\orbit(x_0)\bigr).}
Therefore,
\Eq{*}{
 \sup_{n,m\in\N}d(T^nx_0,T^mx_0)\le\phi\bigl(\diam\orbit(x_0)\bigr)<\diam\orbit(x_0)}
yielding
\Eq{*}{
 \diam\orbit(x_0)=\sup_{n\in\N}d(x_0,T^nx_0).}
On the other hand, $T^nx_0\to x_0$. This property, together with the above one, implies that there exists some
$n_0\in\N$ such that
\Eq{*}{
 \diam\orbit(x_0)=\max\{d(x_0,T^kx_0)\mid k=1,\ldots,n_0\}.}
Let $k\in\{1,\ldots,n_0\}$ be the index via the diameter is represented. Then, for all $n\in\N$, we have
\Eq{*}{
 d(x_0,T^kx_0) & \le d(x_0,T^{n+k}x_0)+d(T^{n+k}x_0,T^kx_0)\\ &
  \le d(x_0,T^{n+k}x_0)+\phi^k\bigl(\diam\orbit(T^nx_0,x_0)\bigr)\\ &
   =d(x_0,T^{n+k}x_0)+\phi^k\bigl(\diam\orbit(x_0)\bigr)\\ &
    \le d(x_0,T^{n+k}x_0)+\phi\bigl(\diam\orbit(x_0)\bigr).}
Passing the limit,
\Eq{*}{
 \diam\orbit(x_0)=d(x_0,T^kx_0)\le\phi\bigl(\diam\orbit(x_0)\bigr)<\diam\orbit(x_0).}
This contradiction completes the proof of the first statement. The second statement is obvious.
\end{proof}

Not claiming completeness, we demonstrate the efficiency of the main result in form of alternative approaches
to known fixed point results. The first one concerns linear quasicontractions of \'Ciri\'c-type \cite{Cir74},
while the second one nonlinear contractions studied by Boyd and Wong \cite{BoyWon69}, by Browder \cite{Brow68},
and by Matkowski \cite{Mat75}. The Banach Fixed Point Theorem is not detailed, since it is an immediate consequence
of both corollaries.

\Corn{If $(X,d)$ is a complete metric space, $q\in]0,1[$ is fixed, and $T\colon X\to X$ satisfies
\Eq{*}{
 d(Tx,Ty)\le q\diam\{x,y,Tx,Ty\}}
for all $x,y\in X$, then $T$ has a unique fixed point. Moreover, the sequence $(T^nx)$ converges to the fixed
point for all $x\in X$.}

\begin{proof}[Hint]
Obviously, $T$ is a strong quasicontraction with comparison function $\phi(t)=qt$. Therefore, in view of the main
result, one should check only the boundedness of orbits. For any $x\in X$, consider its $n$-orbit
\Eq{*}{
 \orbit_n(x)=\{T^kx\mid k=0,1,\ldots,n\}.
}
Then,
\Eq{*}{
 \diam\orbit_n(x)\le d(x,Tx)+\diam\orbit_{n-1}(Tx)
  \le d(x,Tx)+q\diam\orbit_n(x).
}
The arranged form of this inequality shows that $n$-orbits are uniformly bounded. That is, the entire orbit is
bounded, as well.   
\end{proof}

\Corn{If $(X,d)$ is a complete metric space, $\phi$ is a comparison function, and $T\colon X\to X$ satisfies
\Eq{*}{
 d(Tx,Ty)\le\phi\bigl(d(x,y)\bigr)
}
for all $x,y\in X$, then $T$ has a unique fixed point. Moreover, the sequence $(T^nx)$ converges to the fixed
point for all $x\in X$.}

\begin{proof}[Hint]
We shall concentrate again on boundedness of orbits. For this, we use the original (and beautiful) idea of domain
invariance: If $T$ makes small perturbation on the center of a ball, then maps the ball into itself. Let $p\in X$
and $r>0$ be fixed. If $q\in B(p,r)$, then
\Eq{*}{
 d(p,Tq)\le d(p,Tp)+d(Tp,Tq)\le d(p,Tp)+\phi\bigl(d(p,q)\bigr)<r
}
provided that $d(p,Tp)<r-\phi(r)$ holds. That is, $Tq\in B(p,r)$. In particular, $d(p,Tp)<r$; substituting $q=Tp$,
we arrive at $T^2p\in B(p,r)$. Applying induction, $T^mp\in B(p,r)$ follows for all $m\in\N$.

On the other hand, due to the auxiliary lemmas, $d(T^nx,T^{n+1}x)\to 0$ for all $x\in X$. The domain invariance
ensures that $T^{n+m}x\in B(T^nx,r)$ for all $m\in\N$ with a sufficiently large index $n$. In other words, our
mapping induces bounded orbits.  
\end{proof}

Unfortunately, a direct and common generalization of the previous corollaries cannot be given for \emph{arbitrary}
strong quasicontractions. This phenomenon is related to our earlier comment concerning strong quasicontractions that
are not weak ones, since they induces unbounded orbits. Among these kind of strong quasicontractions, there exist
fixed point free ones. However, under an extra assumption on the comparison function, this problem can be avoided,
and the next result is obtained. Note that is also covers the case of Banach. For details of the proof, consult
\cite{Bes15a}.

\Corn{If $(X,d)$ is a complete metric space, $\phi$ is a comparison function such that $\phi^n(t)\le c_nt$ holds with
some convergent series $\sum c_n$, then any strong $\phi$-quasicontraction has a unique fixed point. Moreover, the
sequence of iterates converges pointwise to this fixed point.}

An other straightforward consequence of the main result is that \emph{any weak or strong quasicontraction of a compact
metric space has a unique fixed point}. This can be considered as a counterpart of the statement on strictly nonexpansive
mappings of compact metric spaces.

Finally, let us mention two possibilities on further research. First, as it is well known, the Banach Fixed Point Theorem
has a particular importance in Fractal Theory (see the paper of Hutchinson \cite{Hut81}). A natural question is, what
kind of new impact has the main result in this field? Second, that relaxing the properties of the embedding space is
intensively investigated in present technical literature. For example, the Matkowski Fixed Point Theorem remains true in
so-called regular semimetric spaces \cite{BesPal16}. It is an open problem, whether the main result can be replaced into
this general context or not?

The main theorem points out that boundedness of orbits has more importance then the classical approaches suggest. We
should say that it is a well-established assumption: provides brief form and effective applications simultaneously.
In our opinion, it gives certainly a deeper understanding of the well-known Contraction Principle.

\end{document}